\documentclass[a4paper, 12pt]{article}
\usepackage[utf8]{inputenc}
\usepackage{a4}
\usepackage[T1]{fontenc}
\usepackage{color}
\usepackage[english]{babel}
\usepackage{extarrows}
\usepackage{amssymb,amsmath,amsthm}
\usepackage{enumerate}
\usepackage{graphicx}
\usepackage{epsfig}
\usepackage{comment}
\usepackage{xspace}
\usepackage{float}
\usepackage{multirow}
\usepackage{mathtools}
\usepackage[normalem]{ulem}
\usepackage{color}
\usepackage{nomencl}
\usepackage{nicefrac}
\usepackage{cite}
\makenomenclature
\newcommand{\R}{\ensuremath{\mathbb{R}}\xspace}

\newcommand{\eps}{\epsilon}

\renewcommand{\epsilon}{\varepsilon}

\newcommand{\palpha}{\partial_\alpha}
\newcommand{\udel}{u_\delta}

\newcommand{\dx}{\, \mathrm{d}x}

\newcommand{\intom}{\intop_{\Omega}}
\newcommand{\nabu}{\nabla u}
\newcommand{\tr}{\mathrm{trace}\,}
\newcommand{\gr}[1]{(\ref{#1})} 
\newcommand{\bc}{\big(} 
\newcommand{\bd}{\big)} 
\newcommand{\leb}{\mathcal{L}^n}
\newcommand{\lamh}{\frac{\lambda}{2}}
\begin{document}

\numberwithin{equation}{section}
\newtheoremstyle{break}{15pt}{15pt}{\itshape}{}{\bfseries}{}{\newline}{}
\theoremstyle{break}
\newtheorem*{Satz*}{Theorem}
\newtheorem*{Rem*}{Remark}
\newtheorem*{Lem*}{Lemma}
\newtheorem{Satz}{Theorem}[section]
\newtheorem{Rem}{Remark}[section]
\newtheorem{Lem}{Lemma}[section]
\newtheorem{Prop}{Proposition}[section]
\newtheorem{Cor}{Corollary}[section]
\theoremstyle{definition}
\newtheorem{Def}{Definition}[section]
\newtheorem{exmp}{Example}[section]
\parindent2ex

\newenvironment{rightcases}
  {\left.\begin{aligned}}
  {\end{aligned}\right\rbrace}

\begin{center}{\Large \bf Convex Regularization of Multi-Channel Images Based on Variants of the TV-Model}
\end{center}

\begin{center}
Martin Fuchs, Jan M\"uller, Christian Tietz, Joachim Weickert
\end{center}

\noindent \\
AMS classification: 49J45, 49Q20, 49N60

\noindent \\
Keywords: variational problems of linear growth, TV-regularization, matrix-valued problems

\begin{abstract}
We discuss existence and regularity results for multi-channel images in the setting of isotropic and anisotropic variants of the TV-model.
\end{abstract}
\begin{section}{Introduction}
 In our note we consider a multi-channel image
\[
 f:\Omega\rightarrow\R^N,\;f=(f^1,...,f^N),\,N\geq 1,
\]
defined on a bounded Lipschitz domain $\Omega\subset\R^n$, $n\geq 1$, and try to denoise $f$ by applying a minimization procedure
\begin{align}\label{1.1}
 I[u]:=\intom F(\nabu)\dx+\frac{\lambda}{2}\intom|u-f|^2\dx \rightarrow\min,
\end{align}
where the minimizer $u:\Omega\rightarrow\R^N$ is sought in a suitable class of mappings $w:\Omega\rightarrow\R^N$ depending on the growth of the prescribed density $F$.
We will mainly concentrate on more regular variants $F(\nabu)$ (being convex 
and of linear growth) of the total variation TV-density $|\nabu|$. 
More precisely, we consider two cases: 
\begin{itemize}
\item The {\textbf{isotropic} case:} 
       $F(\nabu)=\varphi\big(\tr(\nabu\nabu^T)\big)=\varphi\left(\sum\limits_{i=1}^n\lambda_i\right).$
\item {The \textbf{anisotropic} case:} 
      $F(\nabu)=\tr\varphi\left(\nabu\nabu^T\right)=\sum\limits_{i=1}^n\varphi(\lambda_i).$
\end{itemize}
Here $\varphi:[0,\infty)\to [0,\infty)$ is a given convex function of linear growth, $\nabla u$ denotes the Jacobian matrix of $u$, $\nabla u^T$ its transpose and $\lambda_1,...,\lambda_n$ are the eigenvalues of the symmetric matrix $\nabla u\nabla u^T$.
The notions of isotropy and anisotropy are motivated by the corresponding
gradient descent evolutions which are diffusion--reaction equations: the 
isotropic setting leads to a diffusion process with a scalar-valued 
diffusivity, while the anisotropic case uses a matrix-valued 
diffusion tensor \cite{WS00b}.

Our goal is to prove existence, regularity and approximation results in both cases for functions $\varphi:[0,\infty)\rightarrow[0,\infty)$ of the principle form
\begin{align}\label{1.2}
 \varphi(s)=\left\{\begin{aligned}
             &\sqrt{\eps^2+s}-\eps,\;\eps>0 \\
               &\hspace{1.3cm}\text{ or }\\
             &\Phi_\mu(\sqrt{s}),\;\mu>1
            \end{aligned}\right\},\;s\geq 0.
\end{align}
Here, 
\[
\Phi_\mu(t):=\intop_0^t\intop_0^s(1+r)^{-\mu}\,dr\,ds,\;p\in \R^{n\times N},\;t\geq 0
\]
 is a standard example for what we call a $\mu$-elliptic density; we refer to  Section 2 and Section 3, respectively, for details including an explanation of the terminology.

The main novelties of this paper can be summarized as follows: for the isotropic case studied in Section 2 we first prove that the relaxed variant of \gr{1.1} admits a unique solution $u$ and in addition we show that the convex-hull-property holds. Second, if the energy density $F$ is $\mu$-elliptic with exponent $\mu <2$, then the solution $u$ is a classical one, i.e. of class $C^{1,\alpha}$. If the anisotropic case is considered (see Sections \ref{Section3} and \ref{Sec4}), then even the proof of the existence of a solution to the relaxed variant of \gr{1.1} is much more elaborate and relies heavily on results of Ball (compare \cite{Ba2} and \cite{Ba1}) on the convexity of functions depending on the singular values of a symmetric matrix. In the same spirit, differentiability results obtained in \cite{Ba3}, enable us to show that the relaxed solution under some assumptions on the data is at least of Sobolev class $W^{1,1}$.

Let us now have a closer look on the history of the problem and its variants 
studied in image analysis. While there has been a long tradition of using 
regularization methods in the context of ill-posed problems \cite{TA77}, 
early quadratic regularization approaches for image analysis problems go 
back to the 1980s \cite{BPT88}. These concepts have been 
generalized to energies with non-quadratic regularizing functions $F$ 
that are either
convex \cite{Sch94} or nonconvex \cite{No90a}. They can be related to the 
nonlinear diffusion filter of Perona and Malik \cite{PM90}. It is fairly
straightforward to extend this diffusion filter to vector-valued
images in the isotropic case \cite{GKKJ92} and to establish corresponding
energies. For matrix-valued data sets, isotropic nonquadratic models
have been pioneered in \cite{TD01a}. Anisotropic regularization approaches
for vector-valued images have been introduced in \cite{WS00b}, and their
matrix-valued counterparts have been considered first in \cite{WB02}. 

The popular TV-regularization approach of Rudin et al. \cite{ROF92} uses 
the total variation seminorm as regularizing function $F(\nabla u)$. 
An early extension of the 
TV-regularizer to color images has been considered in \cite{Ch94}, and
numerous variations of this idea using different channel couplings have 
been proposed within the last two 
decades; see e.g. \cite{DMSC16} and the references therein. A TV-regularization approach for matrix-valued images goes back to \cite{CLLS07}.
In \cite{BBDN06} an anisotropic but rotationally invariant extension of the 
TV-regularizer has been introduced. For further references and a review 
on the large body of work on TV-regularization in image analysis we refer 
to \cite{BO13}.

\end{section}

\begin{section}{Isotropic Regularization}\label{Section2}
 In this section we discuss the following version of the variational problem (\ref{1.1})
\begin{align}\label{2.1}
 J[u]:=\intom\varphi\big(\tr(\nabu\nabu^T)\big)\dx+\frac{\lambda}{2}\intom|u-f|^2\dx\rightarrow\min\text{ in }W^{1,1}(\Omega,\R^N)
\end{align}
with a given function $f:\Omega\rightarrow\R^N$ for which we require
\begin{align}\label{2.2}
 f\in L^2(\Omega,\R^N).
\end{align}
We recall that $\Omega$ is a bounded Lipschitz domain in $\R^n$ and that $\lambda$ denotes some positive number. In what follows, $|\,\cdot\,|$ is the Euclidean norm of vectors and matrices, in particular we have $|\nabu|=\tr\big(\nabu\nabu^T\big)^{\nicefrac{1}{2}}$ for $\nabu=\big(\palpha u^i\big)_{1\leq \alpha\leq n}^{1\leq i\leq N }$. Hence we can write $J[u]$ as
\[
J[u]=\intom\psi\big(|\nabu|\big)\dx+\frac{\lambda}{2}\intom|u-f|^2\dx
\]
with $\psi(s):=\varphi(s^2)$. On the data $f$ we can even impose an extra side condition like
\begin{align}\label{2.3}
 f(x)\in K
\end{align}
for a closed convex subset of $\R^N$, e.g. we can study the case ($N=m^2$, $\R^{m\times m}:=$space of $(m\times m)$-matrices)
\begin{align*}
 K=\mathbb{S}^m:=\Big\{A=(a_{ij})_{1\leq i,j\leq m}\in \R^{m\times m}:a_{ij}=a_{ji},\,i,j=1,...,m, \\ \sum_{i,j=1}^ma_{ij}\xi_i\xi_j\geq\alpha|\xi|^2\text{ for all }\xi\in\R^m\Big\},
\end{align*}
where $\alpha\geq 0$ is fixed. Thus $K$ consists of all symmetric $(m\times m)$-matrices $A$ being $\alpha$-positive (semi-)definite. Concerning the density $\psi:[0,\infty)\rightarrow[0,\infty)$ our assumptions are as follows (and of course partially can be weakened, compare section 3):
\begin{align}
 &\psi\in C^2\big([0,\infty)\big), \psi(0)=0,\label{2.4}\\
 &\psi'(y)\leq \nu_1,\label{2.5}\\
 &\psi(y)\geq\nu_2y-\nu_3,\label{2.6}\\
 &\psi''(y)>0\label{2.7}
\end{align}
for all $y\in [0,\infty)$ and with constants $\nu_1,\nu_2>0$, $\nu_3\in \R$.
Thus
\begin{align}\label{2.8}
 F(p):=\psi\big(|p|\big), p\in \R^{n\times N},
\end{align}
is a strictly convex energy density of linear growth including examples like ($\eps>0$)
\begin{align}\label{2.9}
 F(p):=\sqrt{\eps^2+|p|^2}-\eps,\;p\in\R^{n\times N},
\end{align}
and ($\mu>1$)
\begin{align}\label{2.10}
 F(p):=\Phi_\mu\big(|p|\big),\;\Phi_\mu(t):=\intop_0^t\intop_0^s(1+r)^{-\mu}\,dr\,ds,\;p\in \R^{n\times N},\;t\geq 0.
\end{align}
Recall that we have the following explicit formulas for the functions $\Phi_\mu$
\begin{align}\label{2.11}
 \left\{\begin{aligned}
  &\Phi_\mu(t)=\frac{1}{\mu-1}t+\frac{1}{\mu-1}\frac{1}{\mu-2}(t+1)^{-\mu+2}-\frac{1}{\mu-1}\frac{1}{\mu-2},\;\mu\neq 2,\\
  &\Phi_2(t)=t-\ln(1+t),\;t\geq 0,
 \end{aligned}\right.
\end{align}
and from (\ref{2.11}) we infer that $\Phi_\mu$ approximates the TV-density in the sense that
\begin{align}\label{2.12}
 \lim_{\mu\rightarrow\infty}(\mu-1)\Phi_\mu\big(|p|\big)=|p|,\;p\in\R^{n\times N}.
\end{align}
As a matter  of fact $-$ under the above assumptions on the data $-$ problem (\ref{2.1}) in general fails to have a solution in the Sobolev space $W^{1,1}(\Omega,\R^N)$ and we therefore pass to the relaxed variant of (\ref{2.1}) formulated in the space $BV(\Omega,\R^N)$ of vector-valued functions with finite total variation (see e.g. \cite{AFP}, \cite{Giu} for a definition and further properties of this space). The relaxed variational problem then reads
\begin{align}\label{2.13}
\begin{split}
 K[w]:=&\intom \psi\big(|\nabla^a w|\big)\dx+\psi'_\infty\cdot|\nabla^sw|(\Omega)+\frac{\lambda}{2}\intom|w-f|^2\dx\\ &\rightarrow\min \text{ in }BV(\Omega,\R^N),\quad\psi'_\infty:=\lim_{y\rightarrow\infty}\psi'(y)\in (0,\infty),
\end{split}
\end{align}
where $\nabla w=\nabla^a w\mathcal{L}^n+\nabla^sw$ is the Lebesgue decomposition of the tensor-valued Radon measure $\nabla w$ in its regular and singular part w.r.t. Lebesgue's measure $\mathcal{L}^n$. For details concerning the relaxation procedure the reader is referred e.g. to \cite{BF1,BF2,BF3,FT,FM,Ti}.  We wish to note that 
\begin{align}\label{2.14}
 J[w]=K[w] \text{ for all } w\in W^{1,1}(\Omega,\R^N),
\end{align}
moreover, by standard embedding theorems (compare \cite{Ad} and \cite{AFP}) the finiteness of $\intom|w-f|^2\dx$ for $w\in W^{1,1}(\Omega,\R^N)$ or $w\in BV(\Omega,\R^N)$ is only guaranteed if $n=2$. Let us now state our first result:
\begin{Satz}\label{Thm2.1}
  Assume that we have (\ref{2.2}), (\ref{2.3}) for the data $f$ with $K\subset\R^N$ closed and convex. Then the minimization problem \eqref{2.13} admits a unique solution $u\in BV(\Omega,\R^N)$. The minimizer respects the side-condition (\ref{2.3}), i.e. we have 
\begin{align}\label{2.15}
 u(x)\in K\text{ for almost all }x\in\Omega.
\end{align}
 Moreover it holds
\begin{align}\label{2.16}
 \inf_{u\in BV(\Omega,\R^N)}K[u]=\inf_{w\in W^{1,1}(\Omega,\R^N)}J[w]
\end{align}
with $J$ defined in (\ref{2.1}).
\end{Satz}

\begin{Rem}
 We emphasize that in \gr{2.13} the unconstrained problem is considered, i.e. we do not impose the condition $w(x)\in K$ a.e. on the comparison functions $w\in BV(\Omega,\R^N)$. It just turns out that the unconstrained minimizer $u$ satisfies a kind of maximum-principle better known as convex-hull-property. If we drop the convexity condition for the set $K$, then \gr{2.15} has to be replaced with $u(x)\in\mathrm{conv}(K)$ a.e. on $\Omega$, where $\mathrm{conv}(K)$ is the convex hull of $K$.
\end{Rem}
 Concerning the regularity of the minimizer we have
\begin{Satz}\label{Thm2.2}
 Under the assumptions and with the notation from Theorem \ref{Thm2.1} we impose the following additional requirements on the data $f$ and $\psi$:
\begin{align}
 \label{2.17} f\in L^\infty(\Omega,\R^N),
\end{align}
\begin{align}\label{2.18}
 \begin{split}
  \left\{\begin{aligned}
   &\nu_4(1+t)^{-\mu}\leq \min\left\{\frac{\psi'(t)}{t},\psi''(t)\right\},\\
   &\max\left\{\frac{\psi'(t)}{t},\psi''(t)\right\}\leq \nu_5\frac{1}{1+t}
  \end{aligned}\right.
 \end{split}
\end{align}
for all $t>0$ with positive constants $\nu_4,\nu_5$ and with exponent $\mu>1$. Then, in the case
\begin{align}
 \label{2.19}\mu<2,
\end{align}
problem \gr{2.1} has a solution in the space $W^{1,1}(\Omega,\R^N)$. Moreover, $u$ has Hölder continuous first derivatives in  the interior of $\Omega$.
\end{Satz}
\begin{Rem}
 \begin{enumerate}[i)]
  \item From \gr{2.18} it follows that the density $F$ introduced in \gr{2.8} is $\mu$-elliptic in the sense of
\begin{align}\label{muell}
\nu_4\big(1+|p|\big)^{-\mu}|q|^2\leq D^2F(p)(q,q)\leq\nu_5\frac{|q|^2}{1+|p|},\;p,q\in\R^{n\times N}. 
\end{align}
We remark that example \gr{2.9} satisfies \gr{muell} with exactly $\mu=3$, whereas $F$ from \gr{2.10} satisfies \gr{muell} precisely with the given value of $\mu$.
\item W.r.t. regularity results the bound on $\mu$ stated in \gr{2.19} is  optimal, since even in the case $n=1=N$ there are counterexamples of singular solutions, if the case $\mu>2$ is considered. We refer the reader to \cite{FMT}.
\end{enumerate}
\end{Rem}
Concerning the proofs we just note that Theorem \ref{Thm2.2} is a direct consequence of the results obtained in \cite{BF2}, \cite{BFT} and \cite{Ti}, whereas the existence part of Theorem \ref{Thm2.1} has been established in a very general framework in \cite{FT}. It therefore remains to justify \gr{2.15} for the unique solution $u$ of problem \gr{2.13}. We need the following elementary observation.
\begin{Lem}\label{Lem2.1}
 Consider a closed convex subset $K$ of $\R^N$ and let $\pi:\R^N\rightarrow K$ denote the nearest-point-projection onto $K$, which means that $y_0:=\pi(y)$ is the unique solution of
\begin{align}
 \label{2.20}|y-y_0|=\inf_{z\in K}|y-z|.
\end{align}
The point $y_0$ is characterized through the variational inequality
\begin{align}
 \label{2.21}(y-y_0)\cdot(v-y_0)\leq 0\quad\forall v\in K.
\end{align}
Moreover, the mapping $\pi$ is non-expansive, which means
\begin{align}
 \label{2.22}|\pi(y)-\pi(y')|\leq |y-y'|\quad \forall y,y'\in\R^N.
\end{align}
\end{Lem}
Note that \gr{2.22} is an immediate consequence of \gr{2.21}.

Now, if $f$ satisfies \gr{2.3}, we obtain from Lemma \ref{Lem2.1}
\begin{align*}
 |\pi(w)-f|=|\pi(w)-\pi(f)|\leq|w-f|
\end{align*}
a.e. on $\Omega$ for any measurable function $w:\Omega\rightarrow\R^N$, thus $\pi(u)=u$ and thereby \gr{2.15} holds for our $BV$-solution of problem \gr{2.13} (recall that we have uniqueness), provided we can show that
\begin{align}
 \label{2.23}
\begin{split}
 \intom\psi\bc|\nabla^a\pi(w)|\bd\dx+\psi'_\infty|\nabla^s\pi(w)|(\Omega)\leq \intom\psi\bc|\nabla^aw|\bd\dx+\psi'_\infty|\nabla^sw|(\Omega)
\end{split}
\end{align}
holds for $w\in BV(\Omega,\R^N)$. Inequality \gr{2.23} can be obtained along the lines of the proof of Theorem 1 in \cite{BF4}, however, since the arguments used in this paper are rather technical, we prefer to give a more direct proof of \gr{2.15}. For $\delta>0$ let $F_\delta(p):=\frac{\delta}{2}|p|^2+F(p)$, $p\in\R^{n\times N}$, with $F$ from \gr{2.8} and consider the unique solution of the problem 
\begin{align}\label{Jdelta}
 J_\delta[w]:=\intom F_\delta(\nabla w)\dx+\frac{\lambda}{2}\intom|w-f|^2\dx\rightarrow \min \text{ in }W^{1,2}(\Omega,\R^N).
\end{align}
From \cite{FT}, (4.14), it follows that $u_\delta$ is a $K$-minimizing sequence converging e.g. in $L^1(\Omega,\R^N)$ and a.e. on $\Omega$ to our $K$-minimizer $u$. As remarked above we deduce from \gr{2.3} and \gr{2.22} the validity of
\begin{align*}
 \intom |\pi(\udel)-f|^2\dx\leq\intom |\udel-f|^2\dx,
\end{align*}
whereas from Lemma B.1 in \cite{BF5} it follows
\begin{align*}
 |\partial_\nu(\pi(\udel))|\leq\mathrm{Lip}(\pi)|\partial_\nu\udel|=|\partial_\nu\udel|
\end{align*}
a.e. on $\Omega$, $\nu=1,...,n$. This yields $(\tilde{u}_\delta:=\pi(\udel))$
\begin{align*}
 |\nabla\tilde{u}_\delta|=\left(\sum_{\nu=1}^n|\partial_\nu\tilde{u}_\delta|^2\right)^{\nicefrac{1}{2}}\leq\left(\sum_{\nu=1}^n|\partial_\nu\udel|^2\right)^{\nicefrac{1}{2}}=|\nabla\udel|
\end{align*}
and the structure of $F_\delta$ finally implies $J_\delta[\tilde{u}_\delta]\leq J_\delta[\udel]$, thus $\tilde{u}_\delta=\udel$ by uniqueness. Recalling the convergence $\udel\rightarrow u$ a.e., $\pi(\udel)=\udel$ implies our claim $u=\pi(u)$.\qed

Coming back to the convergence property of the functions $(\mu-1)\Phi_\mu$ stated in formula \gr{2.12} we have the following approximation property of the regularized problems towards the TV-case.

\begin{Satz}\label{Thm2.3}
 Let $\Psi:=(\mu-1)\Phi_\mu$ with $\Phi_\mu$ from \gr{2.10} and let $u_\mu\in BV(\Omega,\R^N)$ denote the unique minimizer of the functional $K$ defined in \gr{2.13} corresponding to this choice of $\Psi$ (note that $\Psi'_\infty=1$), compare with Theorem \ref{Thm2.1}. Then it holds
\begin{align}\label{2.24}
 \|u_\mu-u\|_{L^p(\Omega,\R^N)}\rightarrow 0\quad\forall p<\frac{n}{n-1}
\end{align}
and
\begin{align}\label{2.25}
 u_\mu\rightharpoondown u\text{ in } L^2(\Omega,\R^N)
\end{align}
as $\mu\rightarrow\infty$, where $u$ is the unique minimizer (``TV-solution``) of the problem
\begin{align}\label{2.26}
 \intom|\nabla u|+\frac{\lambda}{2}\intom|u-f|^2\dx\rightarrow\min\text{ in } BV(\Omega,\R^N).
\end{align}
\end{Satz}
\begin{Rem}\label{Rem2.3}
 Clearly a version of Theorem \ref{Thm2.3} also holds for the choice $\Psi(s)=\sqrt{\eps^2+s^2}-\eps$, $\eps>0$, with corresponding solutions $u_\eps$ for which we have the convergences \gr{2.24} and \gr{2.25} as $\eps\downarrow 0$.
\end{Rem}
\begin{Rem}
 Adopting the ideas presented after formula (3.17) in \cite{BF1} it might be possible to improve the convergences \gr{2.24}, \gr{2.25} towards
\begin{align*}
 \lim_{\mu\rightarrow\infty}\|u_\mu-u\|_{L^2(\Omega,\R^N)}=0.
\end{align*}
\end{Rem}

\noindent\textbf{Proof of Theorem \ref{Thm2.3}:} It holds (see formula \gr{2.11})
\begin{align*}
 K[w]=&\intom |\nabla^aw|\dx-\frac{1}{\mu-2}\intom\bc1+|\nabla^aw|\bd^{-\mu+2}\dx-\frac{1}{\mu-2}\leb(\Omega)\\
 &+|\nabla^sw|(\Omega)+\frac{\lambda}{2}\intom|f-w|^2\dx,\; w\in BV(\Omega,\R^N),\;\mu>2,
\end{align*}
and from $K[u_\mu]\leq K[0]$ we directly infer
\begin{align}\label{2.27}
 \sup_{\mu}\left\{\intom|\nabla^au_\mu|\dx+|\nabla^su_\mu|(\Omega)+\intom|u_\mu-f|^2\dx\right\}<\infty,
\end{align}
where we have used that
\begin{align*}
 \frac{1}{\mu-2}\intom\bc1+|\nabla^au_\mu|\bd^{-\mu+2}\dx\rightarrow0\text{ as }\mu\rightarrow\infty.
\end{align*}
Clearly (quoting $BV$-compactness) we can deduce from \gr{2.27} the existence of $\overline{u}\in BV(\Omega,\R^N)$ such that (at least for a subsequence)
\begin{align}\label{2.28}
 \begin{split}
  \left\{\begin{aligned}
   &\|u_\mu-\overline{u}\|_{L^p(\Omega,\R^N)}\rightarrow 0,\;p<\frac{n}{n-1},\\
   & u_\mu\rightharpoondown\overline{u}\text{ in }L^2(\Omega,\R^N)\text{ and }\\
   & u_\mu\rightarrow \overline{u}\text{ a.e.}
  \end{aligned}\right.
\end{split}
\end{align}
holds as $\mu\rightarrow\infty$. By lower semi-continuity of the total variation and by using Fatou's lemma or quoting $u_\mu\rightharpoondown \overline{u}$ in $L^2$ we find
\begin{align*}
 \intom|\nabla\overline{u}|+\frac{\lambda}{2}\intom|\overline{u}-f|^2\dx&\leq \liminf_{\mu\rightarrow\infty}\left(\intom|\nabla u_\mu|+\frac{\lambda}{2}\intom|u_\mu-f|^2\dx\right)\\
 &=\liminf_{\mu\rightarrow\infty}K[u_\mu]\leq \liminf_{\mu\rightarrow\infty}K[u]\\
 &=\intom|\nabla u|+\frac{\lambda}{2}\intom|u-f|^2\dx,
\end{align*}
where we have used the $K$-minimality of the $u_\mu$. Thus $\overline{u}$ is a TV-minimizer, hence $u=\overline{u}$ by the unique solvability of \gr{2.26} and \gr{2.28} is true not only for a subsequence which proves \gr{2.24} and \gr{2.25}. \qed
\begin{Rem}
 We leave it as an exercise to the reader to show that the statements of Theorems \ref{Thm2.1} and \ref{Thm2.3} remain valid if the quantity $\frac{\lambda}{2}\intom|u-f|^2\dx$ is replaced by $\frac{\lambda}{2}\intom\omega(u-f)\dx$ with $\omega:\R^N\rightarrow[0,\infty)$ being strictly convex, e.g. we may choose
\begin{align*}
 \omega(y):=\sqrt{\eps^2+|y|^2}-\eps,\;\eps>0,\;y\in\R^N,
\end{align*}
or $\omega(y):=|y|^p$ with exponent $p>1$. Of course \gr{2.25} then has to be replaced with $u_\mu\rightharpoondown u$ in $L^{n/(n-1)}(\Omega,\R^N)$ in the first case and $u_\mu\rightharpoondown u$ in $L^q(\Omega,\R^N)$ in the second case, where $q:=\max\{p,\nicefrac{n}{n-1}\}$.
\end{Rem}

\begin{Rem}
 If for a given set of data $f$ it is desirable to have smoothness of the regularizer $u$ on a subset $\Omega'$ of $\Omega$, whereas on the complement of $\Omega'$ non-smoothness of $u$ seems to be natural, then such a behavior can be generated by considering non-autonomous densities of the form
 \begin{align*}
  F(x,\nabu)=\eta(x)\Phi_\mu\big(|\nabu|\big)+\big(1-\eta(x)\big)\Phi_\nu\big(|\nabu|\big)
 \end{align*}
with $\mu\in(1,2)$ and $\nu\in(2,\infty)$ large. Here $\eta$ is a smooth function on $\Omega$ such that $0\leq \eta\leq 1$ and with the property $\eta=1$ on $\Omega'$. For details we refer to the paper \cite{BFW2}. 
\end{Rem}

\begin{Rem}
 We note that our discussion can  easily  be extended to isotropic models of super-linear growth. To be precise we consider the problem (compare \gr{2.1})
\begin{align}
 \label{2.29} \intom\Phi_\mu\bc|\nabla u|\bd\dx+\frac{\lambda}{2}\intom|u-f|^2\dx\rightarrow\min
\end{align}
but now with the choice $\mu\leq 1$, where in case $\mu=1$ the correct class for \gr{2.29} is the Orlicz-Sobolev space $W^{1,h}(\Omega,\R^N)$ generated by the function $h(t):=t\ln(1+t)$, $t\geq 0$, (compare \cite{Ad}) and for values $\mu<1$ problem \gr{2.29} is well posed in the Sobolev space $W^{1,p}(\Omega,\R^N)$, $p:=2-\mu>1$. In both cases \gr{2.29} admits a unique solution $u$ satisfying $u(x)\in K$, if $f$ has this property (with $K\subset\R^N$ closed and convex), moreover, it holds $u\in C^{1,\alpha}(\Omega,\R^N)$ for any $\alpha\in (0,1)$. Some details and further references concerning the superlinear case are presented in \cite{BFW2}.
\end{Rem}
\end{section}

\begin{section}{Anisotropic regularization}\label{Section3}
 We start with some preliminaries concerning the definition of the densities $F$ we now have in mind where for notational simplicity we consider the quadratic case for which $n=N$. The general situation is briefly discussed in Remark \ref{Rem3.1}. For matrices $p\in\R^{n\times n}$ let
\begin{align}\label{3.1}
 J(p):=pp^T\;\big((p^T)_{ij}=p_{ji}\big),
\end{align}
and observe that $J(p)$ is symmetric and positive semidefinite with eigenvalues $0\leq\sigma_1(p)\leq...\leq\sigma_n(p)$. We introduce the numbers
\begin{align}\label{3.2}
 \lambda_i(p):=\sqrt{\sigma_i(p)}
\end{align}
which correspond to the eigenvalues of $\sqrt{J(p)}$ and are known as the \textit{singular values} of the matrix $p$. The following observation of Ball (see Theorem 6.1 in \cite{Ba1}) and compare \cite{Ba2}, Theorem 5.1 on p. 363 for a complete proof in any dimension $n$) is of crucial importance
\begin{Lem}
 \label{Lem3.1}
 Consider a function $\rho:[0,\infty)\rightarrow[0,\infty)$ which is convex and increasing. Then the mapping
\begin{align}
 \label{3.3} F:\R^{n\times n}\rightarrow\R,\;p\mapsto\tr\rho\big(\sqrt{J(p)}\big):=\sum_{i=1}^n\rho\bc\lambda_i(p)\bd,\;p\in\R^{n\times n},
\end{align}
is a convex function on the space $\R^{n\times n}$.
\end{Lem}

\begin{Rem}
\label{Rem3.1} For the sake of notational simplicity, we have restricted ourselves to the case of quadratic matrices. However, we would like to indicate how our results can be adapted to the general case of $n\times N$ matrices with $N\neq n$ with the help of of Lemma \ref{Lem3.1}. 
\begin{enumerate}[i)]
\item First we assume $N<n$. Let $p\in\R^{n\times N}$ and $J(p):=pp^T\in\R^{n\times n}$. As before, we denote the eigenvalues of $\sqrt{J(p)}$ by $\lambda_1(p),...,\lambda_n(p)$ and now define $F:\R^{n\times N}\rightarrow \R$ through the formula
\begin{align*}
F:\R^{n\times N}\rightarrow \R,\;p\mapsto\sum_{i=1}^n\rho\big(\lambda_i(p)\big)
\end{align*}
where $\rho:[0,\infty)\rightarrow [0,\infty)$ is as in Lemma \ref{Lem3.1}. Then we define $\widetilde{F}:\R^{n\times n}\rightarrow\R$ according to \gr{3.3}. Now consider the linear embedding $\mathcal{E}:\R^{n\times N}\rightarrow\R^{n\times n}$, which acts on an $(n\times N)$-matrix $p$ by adding $(n-N)$ zero-columns. Then we observe $pp^T=\mathcal{E}(p)\mathcal{E}(p)^T$ for $p\in\R^{n\times N}$ and the convexity follows from the formula  $F(p)=\widetilde{F}\big(\mathcal{E}(p)\big)$ and the convexity of $\widetilde{F}$.

\item The case $N>n$ can be treated in the same manner:  let now $\mathcal{E}:\R^{n\times N}\rightarrow\R^{N\times N}$ denote the embedding which adds $N-n$ zero-rows to a matrix $p\in\R^{n\times N}$, define $F:\R^{n\times N}\rightarrow\R$ as above and $\widetilde{F}:\R^{N\times N}\rightarrow\R$ according to \gr{3.3} (with ``$n$'' replaced by ``$N$''). Then $\mathcal{E}(p)\mathcal{E}(p)^T=pp^T\oplus \mathbf{0}$, where $\mathbf{0}$ denotes the $(N-n)\times(N-n)$-zero matrix and $\widetilde{F}(\mathcal{E}(p))=F(p)+(N-n)\rho(0)$ is a convex function by Lemma \ref{Lem3.1}, and hence so is $F$. 

\item Since the linear map $\mathcal{E}$ is smooth in both cases, we can apply this strategy  to extend our results concerning differentiability in Section \ref{Sec4} to the non-quadratic case.
\end{enumerate}
\end{Rem}
\begin{Rem}
Note, that the general version of Lemma \ref{Lem3.1} as it is found in \cite{Ba2} states, that if $\varphi:\R^{n}\rightarrow\R$ is symmetric and convex, then $\Phi:\R^{n\times n}\rightarrow\R$, $p\mapsto \varphi(\lambda_1(p),...,\lambda_n(p))$ is also convex. The necessity of symmetry is the reason, why we have to apply the same function $\rho$ to each of the eigenvalues $\lambda_i$  in \gr{3.3}. 
\end{Rem}
\begin{Def}[anisotropic energy densities of linear growth]\label{Def3.1}
 Let $\psi:[0,\infty)\rightarrow[0,\infty)$ denote an increasing and convex function satisfying in addition
\begin{align}
 \label{3.4} c_1t-c_2\leq\psi(t)\leq c_3t+c_4
\end{align}
with constants $c_1,c_3>0$, $c_2,c_4\in\R$. Then the mapping (recall \gr{3.1}-\gr{3.3})
\begin{align}
 \label{3.5}
 F_\psi:\R^{n\times n}\rightarrow\R,\;F_\psi:=\tr\psi\big(\sqrt{J}\big),
\end{align}
is termed the anisotropic energy density of linear growth generated by $\psi$.
\end{Def}
This terminology is justified by 
\begin{Lem}
 \label{Lem3.2} In the notation of Definition \ref{Def3.1} the convex function $F_\psi:\R^{n\times n}\rightarrow[0,\infty)$ satisfies
\begin{align}
 \label{3.6} c_1^*|p|-c_2^*\leq F_\psi(p)\leq c_3^*|p|+c_4^*,\; p\in \R^{n\times n},
\end{align}
with constants $c_1^*,c_3^*>0$, $c_2^*,c_4^*\in\R$, $|p|$ denoting the Euclidean (=Frobenius) norm of the matrix $p$.
\end{Lem}
 
\noindent\textbf{Proof of Lemma \ref{Lem3.2}:} From \gr{3.5} together with \gr{3.4} it follows
\begin{align*}
 \sum_{i=1}^n\bc c_1\lambda_i(p)-c_2\bd\leq \sum_{i=1}^n\psi\bc\lambda_i(p)\bd\leq\sum_{i=1}^n\bc c_3\lambda_i(p)+c_4\bd,
\end{align*}
hence
\begin{align*}
 c_1\left(\sum_{i=1}^n\lambda_i(p)\right)-nc_2\leq F_\psi(p)\leq c_3\left(\sum_{i=1}^n\lambda_i(p)\right)+nc_4.
\end{align*}
We further observe 
\begin{align*}
 \sum_{i=1}^n\lambda_i^2(p)=\sum_{i=1}^n\sigma_i(p)=\tr(pp^T)=|p|^2,
\end{align*}
which means
\begin{align*}
 c_5\sum_{i=1}^n\lambda_i(p)\leq |p|\leq c_6\sum_{i=1}^n\lambda_i(p)
\end{align*}
with positive numbers $c_5$, $c_6$. This immediately implies \gr{3.6}.\qed
\begin{exmp}[anisotropic TV-density]
Letting $\psi(t):=t$, $t\geq 0$, in formula \gr{3.5} we obtain 
\begin{align}
 \label{3.7} F_{TV}(p)=\sum_{i=1}^n\lambda_i(p),\;p\in\R^{n\times n}.
\end{align}
Note that the isotropic TV-density is just the quantity $|p|=\big( \sum_{i=1}^n\lambda_i(p)^2 \big)^{\nicefrac{1}{2}}$.
 \label{Ex3.1}
\end{exmp}
\begin{exmp}[regularized TV-densities]
 \label{Ex3.2}
 For $\mu>1$ we let $\psi(t):=\Phi_\mu(t)$, $t\geq 0$, with $\Phi_\mu$ from \gr{2.10} and define 
\begin{align}
 \label{3.8} F_\mu:=(\mu-1)\tr\Phi_\mu\bc\sqrt{J}\bd.
\end{align}
With a slight abuse of notation we can also consider
\begin{align}
 \label{3.9} F_\eps:=\tr\sqrt{\eps^2+J},\;\eps>0
\end{align}
which means that $\psi_\eps(t):=\sqrt{\eps^2+t^2}$ in formula \gr{3.5}.
\end{exmp}
Let us now discuss variational problems in the anisotropic linear growth setting: as usual we consider data 
\begin{align}
 \label{3.10} f\in L^2(\Omega,\R^n)
\end{align}
 for a bounded Lipschitz domain $\Omega\subset\R^n$. For $u:\Omega\rightarrow\R^n$ we let 
\begin{align*}
 \nabla u:=\bc\nabla u^1...\nabla u^n\bd=
\begin{pmatrix}
\partial_1 u^1 & \hdots &\partial_1 u^n\\
\vdots & &\vdots\\ 
\partial_n u^1&\hdots &\partial_n u^n
\end{pmatrix}
\end{align*}
whenever this $(n\times n)$-matrix is defined (in a weak sense). We have (compare \gr{3.1})
\begin{align*}
 J(\nabu)=\nabu\nabu^T=\bc\partial_iu\cdot\partial_ju\bd_{1\leq i,j\leq n},
\end{align*}
''$\cdot$`` denoting the scalar product in $\R^n$, and by Lemma \ref{Lem3.2} the variational problem
\begin{align}
 \label{3.11}J_\psi[u]:=\intom F_\psi(\nabu)\dx+\lamh\intom|u-f|^2\dx\rightarrow\min
\end{align}
is well defined on the Sobolev space $W^{1,1}(\Omega,\R^N)$ for any function $\psi$ as in Definition \ref{Def3.1} and for arbitrary choice of $\lambda>0$. As explained in Section 2 we have to pass to the relaxed version of \gr{3.11} which reads ($\frac{\nabla^sw}{|\nabla^sw|}$ denoting the density of the measure $\nabla^sw$ with respect to the measure $|\nabla^s w|$)
\begin{align}
 \label{3.12}K_\psi[w]:=\intom F_\psi(\nabla^aw)\dx+\intom F_\psi^\infty\left(\frac{\nabla^sw}{|\nabla^sw|}\right)\,\mathrm{d}|\nabla^sw|+\lamh\intop|w-f|^2\dx\rightarrow\min
\end{align}
 in $BV(\Omega,\R^n)$. Here our notation is introduced after \gr{2.13}, and we refer the reader to Theorem 5.47 (and the subsequent remarks) in \cite{AFP}, in particular,
\begin{align*}
 F_\psi^\infty(p):=\lim_{t\rightarrow\infty}\frac{F_\psi(tp)}{t},\;p\in\R^{n\times n},
\end{align*}
is the recession function of $F_\psi$, which here takes the form (compare \gr{3.7})
\begin{align}\label{3.13}
 F_\psi^\infty(p)=\lim_{t\rightarrow\infty}\frac{\psi(t)}{t}\sum_{i=1}^n\lambda_i(p)=\lim_{t\rightarrow\infty}\frac{\psi(t)}{t}F_{TV}(p).
\end{align}
Noting that $\lambda_i\left( \frac{p}{|p|} \right)=\frac{1}{|p|}\lambda_i(p)$, we may therefore write for $w\in BV(\Omega,\R^n)$
\begin{align}\label{3.14}
\begin{split}
 &\intom F_\psi^\infty\left(\frac{\nabla^sw}{|\nabla^sw|}\right)\mathrm{d}|\nabla^sw|=\lim_{t\rightarrow\infty}\frac{\psi(t)}{t}\intom \sum_{i=1}^n\lambda_i\left( \frac{\nabla^sw}{|\nabla^sw|} \right)\,\mathrm{d}|\nabla^sw|\\
 &=\lim_{t\rightarrow\infty}\frac{\psi(t)}{t}\intom \left(\sum_{i=1}^n\lambda_i\right)(\nabla^sw)=\lim_{t\rightarrow\infty}\frac{\psi(t)}{t}\left(\sum_{i=1}^n\lambda_i\right)(\nabla^sw)(\Omega),
\end{split}
\end{align}
where in the last line we apply the convex function $F_{TV}$ (compare \gr{3.7}) to the matrix-valued measure $\nabla^sw$ in the sense of \cite{DT} and calculate the total mass of the resulting nonnegative measure. For the particular case $\psi(t)=t$ the functional \gr{3.12} reduces to
\begin{align}
\begin{split}
 \label{3.15} &K_{TV}[w]=\intom\sum_{i=1}^n\lambda_i(\nabla^aw)\dx+\intom\left(\sum_{i=1}^n\lambda_i\right)(\nabla^sw)+\lamh\intom|w-f|^2\dx,\\ &w\in BV(\Omega,\R^n).
\end{split}
\end{align}
We further like to remark that in formulas \gr{3.13} and \gr{3.14} the quantity $\lim\limits_{t\rightarrow\infty}\frac{\psi(t)}{t}$ can be replaced by (compare \gr{2.13})
\begin{align*}
 \psi'_\infty:=\lim_{t\rightarrow\infty}\psi'(t),
\end{align*}
provided $\psi$ satisfies \gr{2.4}-\gr{2.7}. After these preparations we can state
\begin{Satz}\label{Thm3.1}
 Let $\psi:[0,\infty)\rightarrow[0,\infty)$ denote a convex and increasing function of linear growth as stated in \gr{3.4} and consider the density
\begin{align*}
 F_\psi(p):=\tr\psi\big(\sqrt{pp^T}\big),\;p\in\R^{n\times n},
\end{align*}
being defined in formulas \gr{3.3} and \gr{3.5}.
\begin{enumerate}[a)]
 \item The variational problem (see \gr{3.12})
\begin{align*}
 K_\psi\rightarrow \min\text{ in }BV(\Omega,\R^n) 
\end{align*}
admits a unique solution $u\in BV(\Omega,\R^n)$. It holds
\begin{align*}
 K[u]=\inf_{v\in W^{1,1}(\Omega,\R^n)}J_\psi[v]
\end{align*}
with $J_\psi$ from \gr{3.11}.
\item Let $\psi:=(\mu-1)\Phi_\mu$, i.e. $F_\psi=F_\mu$ with $F_\mu$ from \gr{3.8}, where $\mu>1$. Consider the corresponding version of \gr{3.12}, i.e.
\begin{align*}
 &K_\mu[w]:=\intom(\mu-1)\sum_{i=1}^n\Phi_\mu\big(\lambda_i(\nabla^aw)\big)\dx+\intom\left( \sum_{i=1}^n\lambda_i \right)(\nabla^sw)+\lamh\intom|w-f|^2\dx\\
 &\rightarrow\min\text{ in }BV(\Omega,\R^n)
\end{align*}
with unique solution $u_\mu$. Then it holds
\begin{align*}
 \|u_\mu-u\|_{L^p(\Omega,\R^n)}\rightarrow 0,\;p<\frac{n}{n-1},
\end{align*}
\begin{align*}
 u_\mu-u\rightharpoondown 0\text{ in } L^2(\Omega,\R^n)\text{ and a.e.}
\end{align*}
as $\mu\rightarrow\infty$, where $u\in BV(\Omega,\R^n)$ is the unique TV-solution, i.e. the unique minimizer of the energy $K_{TV}$ defined in \gr{3.15}.
\item For $\eps>0$ let $\psi(t):=\psi_\eps(t):=\sqrt{\eps^2+t^2}$, $t\geq 0$, in \gr{3.12}, i.e. we look at the problem
\begin{align*}
 &K_\eps[w]:=\intom\sum_{i=1}^n\sqrt{\eps^2+\lambda_i(\nabla^aw)^2}\dx+\intom\left( \sum_{i=1}^n\lambda_i \right)(\nabla^sw)+\lamh\intom|w-f|^2\dx\\
 &\rightarrow\min\text{ in }BV(\Omega,\R^n)
\end{align*}
with corresponding solution $u_\eps$. Then we have 
\begin{align*}
 \|u_\eps-u\|_{L^p(\Omega,\R^n)}\rightarrow 0,\;p<\frac{n}{n-1},
\end{align*}
\begin{align*}
 u_\eps-u\rightharpoondown 0\text{ in } L^2(\Omega,\R^n)\text{ and a.e.}
\end{align*}
as $\eps\rightarrow 0$, where $u$ is the solution of (see \gr{3.15})
\begin{align*}
 K_{TV}\rightarrow\min\text{ in }BV(\Omega,\R^n). 
\end{align*}
\end{enumerate}
\end{Satz}

\noindent\textbf{Proof of Theorem \ref{Thm3.1}:} a) Let $u_k$ denote a $K_\psi$-minimizing sequence from $BV(\Omega,\R^n)$. Lemma \ref{Lem3.2} (compare inequality \gr{3.6}) in combination with the definition of $K_\psi$ then yields
\begin{align*}
 \sup_{k}|\nabla u_k|(\Omega),\;\sup_{k}\|u_k\|_{L^2(\Omega,\R^n)}<\infty,
\end{align*}
hence, quoting $BV$-compactness, $u_k\rightarrow:\overline{u}$ in $L^1(\Omega,\R^n)$ for some $\overline{u}\in BV(\Omega,\R^n)$ and a subsequence of $u_k$. Moreover, we may assume that $u_k\rightarrow\overline{u}$ a.e. on $\Omega$ and 
\begin{align*}
 \intom|\overline{u}-f|^2\dx\leq\liminf_{k\rightarrow\infty}\intom|u_k-f|^2\dx
\end{align*}
follows from Fatou's Lemma (or from $u_k\rightharpoondown\overline{u}$ in $L^2(\Omega)$). According to Theorem 5.47 in \cite{AFP} and the remarks stated after this theorem the functional
\begin{align*}
 w\mapsto\intom F_\psi(\nabla^aw)\dx+\intom F_\psi^\infty\left( \frac{\nabla^sw}{|\nabla^sw|} \right)\,\mathrm{d}|\nabla^sw|
\end{align*}
is lower semi-continuous with respect to $L^1(\Omega,\R^n)$-convergence. Here we make essential use of Ball's convexity result Lemma 3.1 implying the 
convexity of $F_\psi$. 
Altogether we have
\begin{align*}
 K_\psi[\overline{u}]\leq \liminf_{k\rightarrow\infty}K_\psi[u_k],
\end{align*}
thus $\overline{u}$ is $K_\psi$-minimizing. Uniqueness of the minimizer is immediate, all other claims follow along the lines of Theorem \ref{Thm2.1}.

For part $b)$ and $c)$ we refer to Theorem \ref{Thm2.3} and Remark \ref{Rem2.3}. \qed
\end{section}

\begin{section}{Differentiable models}\label{Sec4}
Concerning the regularity properties of the minimizer $u\in BV(\Omega,\R^n)$ from Theorems \ref{Thm2.1} and \ref{Thm3.1}, it is desirable to consider energy densities $F$ which are sufficiently smooth. Namely we would like to have $F\in C^2(\R^{n\times n})$. To this end, we consider a slight modification of the function $F$ from \gr{3.3} by setting
 \begin{align}\label{5.2*}
 F^*(p):=\sum_{i=1}^n\psi\left(\sqrt[4]{\eps^2+\sigma_i(p)^2}\right)
\end{align}
for some $\eps>0$, with $\sigma_i$ as usual denoting the eigenvalues of $pp^T$ and $\psi:\R\rightarrow[0,\infty)$ is a convex and increasing function which satisfies (cf. \gr{2.4}-\gr{2.7}) 
\begin{align}\label{5.1}
\left\{\begin{aligned}
 &\psi\in C^2(\R),\\
 &\psi(-y)=\psi(y),\;\psi(0)=0,\\
 &|\psi'(y)|\leq\nu_1,\\
 &\psi(y)\geq\nu_2|y|-\nu_3,\;\psi''(y)> 0,\\
\end{aligned}\right.
\end{align}
with $\nu_1,\nu_2,\nu_4>0$, $\nu_3\in \R$. As for the map $p\mapsto (\sigma_1(p),...,\sigma_n(p))$, which is not immediately seen to be differentiable, we can once more benefit from a result by John Ball in \cite{Ba3}  which gives us the desired smoothness. Precisely we have
\begin{Satz}\label{Thm5.1}
 The density $F^*$ being defined in \gr{5.2*} is convex and $C^2$ on $\R^{n\times n}$.
\end{Satz}
\begin{Rem}
As we have already mentioned in Remark \ref{Rem3.1} iii), the above result can easily be adjusted to the non-quadratic case $f:\Omega\rightarrow\R^N$ for $N\neq n$. 
\end{Rem}
\noindent\textbf{Proof of Theorem \ref{Thm5.1}:}

\noindent  With the notation from \gr{3.3} and \gr{3.5} we have
\[
 F^*(p)=\tr\widetilde{\psi}\bc\sqrt{J}\bd,\;J=J(p)=pp^T,
\]
i.e. $F^*(p)=\sum_{i=1}^n\widetilde{\psi}\bc\lambda_i(p)\bd$, if we set
\[
 \widetilde{\psi}(t):=\psi\left( \sqrt[4]{\eps^2+t^4} \right).
\]
Since $\widetilde{\psi}$ fulfills the requirements imposed on $\rho$ of Lemma \ref{Lem3.1}, the convexity of $F^*$ follows. We use the notation from \cite{Ba3}, Section 5. Let 
\begin{align}\label{defh}
\begin{split}
\begin{rightcases}
 &E:=\mathbb{S}^n\quad\text{(symmetric $(n\times n)$-matrices)},\\
 &\Gamma_E:=\{\text{diagonal matrices in }E\},\\
 &v_i(A):=\text{i-th eigenvalue of }A\in E,\\
 &H:\Gamma_E\cong \R^n\ni(t_1,...,t_n)\mapsto\sum_{i=1}^n\psi\left(\sqrt[4]{\eps^2+t_i^2}\right).
\end{rightcases}
\end{split}
\end{align}
Obviously, $H\in C^2(\R^n)$. But then, Theorem 5.5 on p. 717 in \cite{Ba3}, implies that also
\[
 h:\mathbb{S}^n\ni A\mapsto H(v_1(A),...,v_n(A))
\]
is of class $C^2$ on $\mathbb{S}^n$ $(\cong \R^{n(n-1)/2})$. Now note that
\[
 F^*(p)=h(pp^T),\; p\in\R^{n\times n}
\]
and since the map $p\mapsto pp^T$ is obviously smooth, this shows $F^*\in C^2(\R^{n\times n})$.\qed

\begin{Rem}
The symmetry of the function $H$ is essential for establishing both convexity and differentiability of our models. In particular we cannot generalize our model to $\sum_{i=1}^n\psi_i\left(\sqrt[4]{\eps+t_i^2}\right)$ with distinct $\psi_i$'s for each $i\in\{1,...,n\}$.
\end{Rem}
\begin{Satz}\label{Thm5.2}
 Let $\psi$ satisfy \gr{5.1} and define $F^*$ according to (\ref{5.2*}).
 \begin{enumerate}[a)]
  \item $F^*$ grows linearly in the sense of inequality \gr{3.6}.
  \item The relaxation of
\begin{align}\label{5.7}
\intom F^*(\nabu)\dx+\lamh \intom|u-f|^2\dx\rightarrow\min\text{ in }W^{1,1}(\Omega,\R^n)
\end{align}
with $f\in L^2(\Omega,\R^n),\;\lambda>0$, is given by
\begin{align}\label{5.8}
\begin{split}
 &\intom F^*(\nabla^a u)\dx+\psi'_\infty\cdot\left( \sum_{i=1}^n\lambda_i \right)(\nabla^su)(\Omega)+\lamh \intom|u-f|^2\dx\rightarrow\min\\
 &\text{ in }BV(\Omega,\R^n)
\end{split}
\end{align}
and is uniquely solvable. (Here we have abbreviated $\psi'_\infty:=\lim\limits_{s\rightarrow\infty}\psi'(s)=\lim\limits_{s\rightarrow\infty}\frac{\psi(s)}{s}$.)
\item Let us in addition  assume that $F^*$ satisfies
\begin{align}\label{D2F*}
\big|D^2F^*(p)\big|\leq \nu_4\frac{1}{1+|p|}
\end{align}
for some constant $\nu_4>0$. Then, if $\Omega'\subset\Omega$ and $f\in W^{1,2}_\mathrm{loc}(\Omega',\R^n)$, we have $u\in W^{1,2}_\mathrm{loc}(\Omega',\R^n)$ for the unique solution $u$ of \gr{5.8}.
 \end{enumerate}
\end{Satz}
\begin{Cor}\label{Cor5.1}
 If the data $f$ are chosen from the space $W^{1,2}_\mathrm{loc}(\Omega,\R^n)$ and $F^*$ satisfies \gr{D2F*}, then \gr{5.7} is solvable in $W^{1,1}(\Omega,\R^n)$.
\end{Cor}
\begin{Rem}
 As usual Theorem \ref{Thm5.2} and Corollary \ref{Cor5.1} extend to the non-quadratic case $f:\Omega\rightarrow\R^N$ for $N\neq n$ via Remark \ref{Rem3.1}.
\end{Rem}

\begin{Rem}
 In $\left( \sum\limits_{i=1}^n\lambda_i \right)(\nabla^su)$ the convex function $p\mapsto\sum\limits_{i=1}^n\lambda_i(p)$ is applied to the matrix-valued measure $\nabla^su$ which yields a positive Radon measure on $\Omega$, whose total mass enters in \gr{5.8}. We refer to the comments after formula \gr{3.14}.
\end{Rem}

\noindent\textbf{Proof of Theorem \ref{Thm5.2}:} Ad a): cf. the proof of Lemma \ref{Lem3.2};

\noindent Ad b): see Theorem \ref{Thm3.1} and note (cf. \gr{3.13}) that
\begin{align*}
 &\bc F^\infty\bd(p):=\lim_{t\rightarrow\infty}\frac{1}{t}F^*(tp)=\psi'_\infty\sum_{i=1}^n\psi\left( \sqrt[4]{\sigma_i^2(p)} \right)=\psi_\infty'\left( \sum_{i=1}^n\lambda_i \right)(p)\\
&\big(\lambda_i(p):=\text{ eigenvalues of }\sqrt{pp^T}=\sqrt{\sigma_i(p)}\big).
\end{align*}
This implies (cf. \gr{3.7})
\[
 \bc F^*\bd^\infty(p)=\psi_\infty'F_{TV}(p).
\]

\noindent Ad c): let w.l.o.g. $f\in W^{1,2}_\mathrm{loc}(\Omega,\R^n)$. In all the following calculations we have to replace $u$ with the sequence of regularizers $\udel$ (cf. \gr{Jdelta} and compare \cite{FT} for more details), however, for notational simplicity we drop the index $\delta$, i.e. $F^*(p)=F^*_\delta(p):=\frac{\delta}{2}|p|^2+ F^*(p)$ and $u=\udel$ is the unique solution of 
\[
 \intom F^*_\delta(\nabla w)\dx+\lamh\intom{|w-f|^2}\dx\rightarrow\min\text{ in }W^{1,2}(\Omega,\R^n).
\]

From the minimality of $u$ along with $F^*\in C^2$ it follows (using summation convention w.r.t. the index $\alpha$)
\begin{align*}
 \intom D^2F^*(\nabu)\bc\palpha\nabu,\nabla(\eta^2\palpha u)\bd\dx=\lambda\intom\palpha(\eta^2\palpha u)\cdot(u-f)\dx,
\end{align*}
where $\eta\in C^\infty_0(\Omega)$, $\text{spt } \eta\subset B_{2R}(x_0)$ with $0\leq \eta\leq 1$ and $\eta\equiv 1$ on $B_R(x_0)$ for some $x_0\in\Omega$ and some radius $R>0$ s.t. $B_{2R}(x_0)\subset\Omega$ . Hence
\begin{align*}
 &\intom D^2F^*(\nabu)\bc\eta\palpha\nabu,\eta\palpha\nabu\bd\dx+\intom D^2F^*(\nabu)\bc\palpha\nabu,\nabla\eta^2\otimes\palpha u\bd\dx\\
 &\quad+\lambda\intom \eta^2|\nabu|^2\dx=\lambda\intom\eta^2\palpha u\cdot \palpha f\dx
\end{align*}
and thus
\begin{align}\label{5.9}
\begin{split}
 &\intom D^2F^*(\nabu)\bc\eta\palpha\nabu,\eta\palpha\nabu\bd\dx+\lambda\intom \eta^2|\nabu|^2\dx\\
&= \lambda\intom\eta^2\palpha u\cdot \palpha f\dx-\intom D^2F^*(\nabu)\bc\palpha\nabu,\nabla\eta^2\otimes\palpha u\bd\dx\\ &=:T_1+T_2.
\end{split}
\end{align}
The integral $T_1$ can be estimated by Young's inequality through
\begin{align*}
 |T_1|\leq  c(\eps,\lambda)\intom \eta^2|\nabla f|^2\dx+\lambda\eps\intom\eta^2|\nabla u|^2\dx\leq c(\eps,\lambda,R)+\lambda\eps\intom\eta^2|\nabla u|^2\dx.
\end{align*}
Choosing $\eps=1/2$ and absorbing terms on the left-hand side of \gr{5.9}, we obtain
\begin{align*}
 &\intom D^2F^*(\nabu)\bc\eta\palpha\nabu,\eta\palpha\nabu\bd\dx+\frac{\lambda}{2}\intom \eta^2|\nabu|^2\dx\leq c(R)+T_2.
\end{align*}
Now, for $T_2$, we apply the Cauchy-Schwarz inequality to the bilinear form $D^2F^*(\nabu)$ observing
\[
D^2F^*(\nabu)\bc\palpha\nabu,\nabla\eta^2\otimes\palpha u\bd=2D^2F^*(\nabu)\bc\eta\palpha\nabu,\nabla\eta\otimes\palpha u\bd
\]
and obtain after an application of Young's inequality the following result:
\begin{align}\label{5.10}
\begin{split}
 |T_2|\leq &c(\eps)\intom D^2F^*(\nabu)\big(\nabla\eta\otimes\palpha u,\nabla\eta\otimes\palpha u\big)\dx\\
&+\eps\intom D^2F^*(\nabu)\big(\eta\palpha\nabu,\eta\palpha\nabu\big)\dx.
\end{split}
\end{align}
Choosing $\eps=1/2$, the second term on the right-hand side of \gr{5.10} can be absorbed in the left-hand side of \gr{5.9}. For estimating the first term on the right-hand side of \gr{5.10}, we need our additional assumption \gr{D2F*} on $D^2F^*$ which yields:
\begin{align}\label{5.11}
\begin{split}
 \intom D^2&F^*(\nabu)\big(\nabla\eta\otimes\palpha u,\nabla\eta\otimes\palpha u\big)\dx\\
 &\leq c\intom D^2F^*(\nabu)\big(\nabu,\nabu\big)\dx\leq c\intom \big(1+|\nabu|\big)\dx
 \end{split}
\end{align}
with a suitable constant $c$ uniformly with respect to the (invisible) parameter $\delta$. Consequently, \gr{5.9} yields a uniform (in $\delta$) bound for $\intom |\nabu|^2\dx$, which concludes the proof of Theorem \ref{Thm5.2}.\qed

\end{section}

\begin{tabular}{l}
Martin Fuchs (fuchs@math.uni-sb.de)\\
Jan M\"uller (jmueller@math.uni-sb.de)\\
Christian Tietz (tietz@math.uni-sb.de)\\
Joachim Weickert (weickert@mia.uni-saarland.de)\\
\vspace{-0.3cm}\\
Saarland University\\
Department of Mathematics\\
P.O. Box 15 11 50\\ 
66041 Saarbr\"ucken\\ 
Germany
\end{tabular}

\end{document}